\documentclass[11pt]{amsart}
\RequirePackage[OT1]{fontenc}
\RequirePackage[colorlinks,citecolor=blue,urlcolor=blue]{hyperref}
\usepackage{color}
\usepackage{graphicx}
\usepackage{enumitem}

\usepackage{bbm}

\def\A{\mathcal A}

\def\R{{\mathbb R}}
\def\C{{\mathbb C}}
\def\N{{\mathbb N}}
\def\Z{{\mathbb Z}}
\def\tr{\;\mathrm{tr}\,}
\def\<{\langle}
\def\>{\rangle}
\def\P{\mathbb P}

\def\E{\mathbb E}

\def\eps{\epsilon}
\def\fif{\varphi}

\def\i{\mathbb I}
\def\0{\underline 0}
\def\ka{\kappa}

\def\G{\mathcal G}
\def\K{\mathcal K}

\def\e{\mathbf e}
\def\mk {\mathcal {MK}}
\def\W{\mathcal W}
\def\1{\underline 1}

\def\x{\mathbf X}
\def\y{\mathbf Y}

\def\v{\mathbf V}
\def\u{\mathbf U}
\def\z{\mathbf Z}

\def\var{\mathbb{V}\mathrm{ar}}

\def \P {{\mathbb P}}

\def \E {{\mathbb E}}
\def \R {{\mathbb R}}

\def\X{\mathbb X}
\def\Y{\mathbb Y}
\def \U {{\mathbb U}}
\def \V {{\mathbb V}}
\def\Z{\mathbb Z}

\def \N {{\mathbb N}}
\def \K {\mathcal K}
\def \k {\kappa}

\newcommand{\bel}{\begin{equation}\label}
\newcommand{\nobel}{\begin{equation}}
\newcommand{\ee}{\end{equation}}
      \newtheorem{theorem}{Theorem}[section]
    \newtheorem{proposition}[theorem]{Proposition}
       \newtheorem{corollary}[theorem]{Corollary}
       \newtheorem{lemma}[theorem]{Lemma}
      \newtheorem{remark}{Remark}[section]
\theoremstyle{definition}
\newtheorem{definition}{Definition}[section]

\begin{document}
%
%
%
%
%

%

\title{REGRESSION CONDITIONS THAT CHARACTERIZE FREE--POISSON AND FREE--KUMMER DISTRUBTIONS}
\author{Agnieszka Piliszek}
\address{Agnieszka Piliszek\\
Wydzia{\l} Matematyki i Nauk Informacyjnych\\
Politechnika Warszawska\\
Koszykowa 75\\
00-662  Warszawa, Poland}
\email{A.Piliszek@mini.pw.edu.pl}
\date{\today}
\maketitle

		\begin{abstract}
	We find the asymptotic spectral distribution of random Kummer matrix. 
Then we formulate and prove a~free analogue of HV independence property, which is known for classical Kummer and Gamma random variables and for Kummer and Wishart matrices.
We also prove a related characterization of free--Kummer and free--Poisson (Marchenko-Pastur) non--commutative random variables.
		\end{abstract}

\section{Introduction}

This paper concerns some connections between classical  and non-commutative probability, 
especially the links between independence and freeness. 
We prove a free-analogue of a classical indpendence property of Kummer and Gamma random variables. 
Similar attempts have been succesfull for the Lukacs' characterization of the Gamma distribution, \cite{BB04, Szp15}, 
as well as for the Matsumoto-Yor characterization of GIG and Gamma distributions, \cite{Szp17}.
Earlier, also Kac--Bernstein characterization of independent Gaussian random variables, \cite{Be41} was proved for semicircle free non-commutative random variables, \cite{Nic96}.
However it is not clear which of independence characterizations  known for commutative random variables, 
would also hold for their non--commutative counterparts (and what would be the counterparts). 
An important property of classical Gaussian random variables, known as Cramer's Theorem, does not hold in a free probability setting, \cite{BV95}.

Let us recall that a random variable $Y$ has the Gamma distribution with parameters $a,c>0$, we write $Y\sim\G(a,\,c)$,
if it has a density function
$$f(y) \propto y^{a-1}\, e^{-cy}\, I_{(0,\infty)}(y).$$
A random variable $X$ has the Kummer distribution with parameters  $a,c>0$, $b\in\R$, we write ${X\sim\K(a,\,b,\,c)}$, 
if it has a density function
$$ f(x) \propto x^{a-1}\,(1+x)^{-(a+b)}\, e^{-cx}\,I_{(0,\infty)}(x).
$$
An interesting property noticed in \cite{HV16} says, that if  
$X\sim\K(a,b,c)$ and $Y\sim\G(a+b, c)$ are independent, 
then random variables 
\bel{uiv} U := \frac{Y}{1+X}\qquad \mathrm{and}\qquad V := X\;(1+U)\ee 
are also independent and $U\sim\K(a+b, -b, c)$, $V\sim\G(a, c)$. 
We call this property the {\bf HV property} refering to the names of its authors, \cite{HV16}. 
In \cite{PW16}, under the assumption that densities of $X$ and $Y$ are locally integrable, the~converse theorem was proved: 
if $X$ and $Y$ are independent and $U$ and $V$, defined in~\eqref{uiv}, are independent, 
then necessarily $X\sim\K(a,b,c)$ and $Y\sim\G(a+b, c)$ for some constants  $a,c>0$ and $b>-a$. 
In \cite{PW18} it was shown that instead of independence of $U$ and $V$, it is enough to assume constant regression conditions: $\E(V|U) =\alpha$ and $\E(V^{-1}|U) =\beta$ 
(in fact there is a whole range of  equivalent conditions, see \cite{PW18}). 
Some of other properties and characterizations of this type (we mean regression conditions) have its counterparts in non-commutative probability:
\begin{itemize}
	\item Kac--Bernstein characterization of Gaussian distribution, \cite{Be41}, and related free characterization of~Wigner law, \cite{Nic96};
	\item Lukacs' Theorem, \cite{Lu55}, which characterizes Gamma distribution through independence of ${X+Y}$ and $X/Y$. 
		In the non-commutative version Marchenko-Pastur (free-Poisson) distribution is characterized by conditions: 
		$$\fif(\X|\X+\Y) = \alpha (\X+\Y) + \alpha_0,\quad\var(\X|\X+\Y) = C(1+a(\X+\Y) + b(\X+\Y)^2)$$
		in \cite{BB04} or by constant conditional moments of order $1$ and $-1$ of $(\X+\Y)^{-1/2}\X(\X+\Y)^{-1/2}$ in~\cite{Szp16}; 
	\item Matsumoto--Yor property and related characterization of GIG and Gamma, \cite{LW00}, 
		also have their counterparts in free probability, cf. \cite{Szp17}.
\end{itemize}
These results, among others, indicate certain link between commutative and non--commutative probability and between independence and freeness. 

In this paper we give one more link of a similar type, which seems to be of interest on its own.
The basic aim of this paper is to find a pair of~distributions $\mu, \nu$ such that if $\X\sim \mu$ and $\Y\sim\nu$ are free, 
self-adjoint random variables from some $C^*$ algebra $\A$, then $\U = (\i+\X)^{-1/2}\,\Y\, (\i+\X)^{-1/2}$ and 
$\V = (\i+\U)^{1/2}\,\X\,(\i+\U)^{1/2}$ are also free. 
One may compare this transformation with \eqref{uiv}, as it is a generalization of \eqref{uiv}.  
It seems justified to use free-Poisson as $\nu$, which played the role of the Gamma distribution in both cases of Lukacs' and Matsumoto-Yor theorems, \cite{Nic96, Szp17}.
A candidate for $\mu$ is considered in Section \ref{sec3}, and free property of $\mu$ and $\nu$ is proven in Section~\ref{sec4}.
Section~\ref{sec5} contains the main result of the paper (Thm. \ref{free_char}): the constant regression characterization of the measures $\mu$ and $\nu$.

\section{Preliminaries}

A broad introduction to free probability can be found  in \cite{HP} or \cite{MS2017}. 
To increase 
accessability of the paper we recall that part of this theory we 
need in order to present our results. 

Let $\mathcal{A}$ be a $\star$-algebra and $\fif:\mathcal{A} \to\C$ be a linear functional such that $\fif(\i) = 1$, where $\i = 1_{\mathcal {A}}$. 
We assume that $\fif$ is faithful, normal, tracial and positive.
We call the pair $(\mathcal{A}, \fif)$ a non--commutative $\star $-probability space and elements of $\mathcal A$ 
are called non-commutative {random variables}
. 

An important case is $\A\subset \mathcal{B}(H)$, where $H$ is a Hilbert space and $\mathcal{B}(H)$ denotes the space of~bounded linear operators from $H$ to $H$. 
If $\A$ is von Neummann algebra, then we say that $(\A,\fif)$ is $W^\star$-probability space.

Let $a_1, \ldots, a_n\in\A$ be bounded. 
The numbers $\fif(a_{i(1)}\cdots a_{i(k)})$, $i(j)\in\{1,\ldots,n\}$ for
$j\in\{1,\ldots,k\}$ are called \textit{moments} and by 
\textit{joint distribution } of $(a_1,\ldots,a_n)$ we mean the collection of all moments. 
For a bounded,  self--adjoint random variable $\X\in\A$ we can define the $\star$-distribution of $\X$ 
as a unique, real, compactly supported probability measure $\mu$. 
This measure $\mu$ is uniquely determined as such that for all $n\in\N$
$$\fif\left(\X^n\right) =\int_{\R} t^n\, d\mu(t).$$ 
If the support of $\mu$ is contained in $(0,\infty)$, then we say that $\X$ is positive.
For a $n$--tuple of self--adjoint random variables $(\X_1, \ldots, \X_n)$  its joint distribution is defined as 
the linear functional $\mu_{(\X_1, \ldots, \X_n)}:
\C\langle x_1,\ldots, x_n\rangle \to \C$
such that for any polynomial $P\in\C\langle x_1,\ldots, x_n\rangle$ in non-commuting variables the following equality holds:
$$ \fif\left(P\left(\X_1, \ldots, \X_n\right)\right) = \mu_{(\X_1, \ldots, \X_n)}(P)
.
$$
 
Let $(\A,\,\fif)$ be a non--commutative probability space and let $I$ be a finite index set.
For each $i\in I$ let $\A_i\subset\A$ be a unital subalgebra. 
The subalgebras $(\A_i)_{i\in I}$ are called {\it free} or {\it freely independent}, if 
$\fif(a_1\cdots a_k) = 0$ whenever the following four conditions hold:	
\begin{enumerate}
	\item $k\geq 2$ is a positive integer, 
	\item $a_j\in \A_{i(j)}$ ($i(j)\in I$) for all $j=1,\ldots, k$,
	\item  $\fif(a_j) = 0$ for all $j=1,\ldots, k$,
	\item neighbouring elements are from different subalgebras, i.e., $i(1)\neq i(2) \neq\ldots\neq i(k)$.
\end{enumerate}
We say that random variables $(\X_i)_{i\in I}$, $\X_i \in \A$ for each  $i\in I$, are {\it free} or {\it freely independent}, if $(\A_i,\, \i)_{i\in I}$ are free, where $\A_i$ is a subalgebra of $\A$ generated by $\X_i$ for each $i\in I$. 
Freeness of non-commutative random variables can be also expressed conveniently in terms of~free--cumulants. 

Let $NC(n)$ denote the set of all non--crossing partitions of the set $\{1,\ldots,n\}$. 
We define {\it free cumulants} $\k_n:\A^n\to\C$, $n\geq 1$ as multi-linear functionals by the recursive moment-cumulant relation
$$\fif(a_1\ldots a_n) = \sum_{\pi\in NC(n)} \k_\pi(a_1,\ldots, a_n),$$
where $\k_\pi$ is a product of cumulants over all blocks of $\pi$ and the arguments are given by the elements corresponding to the respective blocks. 
For example, if $\pi= \{\{1,4\}, \, \{2,3\}\}$ then $\k_\pi(a,b,c,d) = \k_2(a,d)\, \k_2(b,c)$. 
So $\k_1(a) = \fif(a)$, $\k_2(a_1,\, a_2) = \fif(a_1\, a_2) - \fif(a_1)\,\fif(a_2)$ and so on.

Random variables $\X$ and $\Y$ are free if and only if $\k_n(a_1,\ldots, a_n) = 0$ whenever: $n\geq 2$, $a_i\in\{\X,\,\Y\}$ for all $i$ and there are at least two indices $i$, $j$ such that $a_i=\X$, $a_j = \Y$, cf. \cite{Speicher94}. 
Throughout following sections, we will use a well known formula, which connects free cumulants and moments:
\bel{wzorBLS}
\fif(\X_1\; \ldots\; \X_n) = 
\sum_{k=1} ^n \sum_{1<i_2<\ldots<i_k\leq n} \k_k(\X_1, \X_{i_1},\ldots, \X_{i_k})\prod_{j=1}^k\fif (\X_{i_j +1}\;\ldots\;\X_{i_{j+1} -1}),
\ee
where $i_1 =1$, $i_{n+1} = n+1$, cf. \cite{BLS}.

\subsection{Analytical tools}
	Let $(\A, \fif)$ be a non-commutative probability space.
	Let $\X\in\A$.
	Let us recall that the $r$--transform of~a~random variable $\X$ is: 
	$$r_\X(z) = \sum_{n=0}^\infty \k_{n+1}(\X)\, z^n,$$
	where $\k_n(\X) = \k_n(\overbrace{\X, \ldots, \X}^{n\;\mathrm{times}})$.
	It is known that if $\X$ has a compact support, then  its $r$--transform is an~analytic function in the neighborhood of $0$ (as a function of complex variable).
	For a self--adjoint, bounded $\X$ from a $\C^\star$-algebra $\A$ with a state $\fif$ and $\star$-distribution $\nu$ its {Cauchy transform} is
	$$G_\X(z)= G_\nu(z) = \fif\left((z-\X)^{-1}\right) = \int_\R(z-t)^{-1}\,d\nu(t)$$
	for $z\in\C\setminus{\R}$.
	Let $\C^+ = \{z\in\C:\, \Im(z) >0\}$ and $\C^- = \{z\in\C:\, \Im(z)<0\}$.
	Then $G_\X$ is an analytic function from $\C^+$ to $\C^-$ (Lemma 3.1.2 in \cite{MS2017}).
	Cauchy transform uniquely determines the distribution and measure $\nu$ can be recovered form $G_\nu$ by Stieltjes inversion formula: 
	\begin{proposition}[Theorem 3.1.6 in \cite{MS2017}]
		\label{Stieltjes}
		If $\nu$ is a probabilistic measure on $\R$ and $G_\nu$ is its Cauchy transform, then for every $a<b$
		\bel{stielt}
		 -\lim_{\eps\to 0}\frac{1}{\pi}\int_a^b \Im G_\nu(x+i\eps)\, dx = \nu\left((a,b)\right) + \frac{1}{2}\,\nu\left(\{a,b\}\right).
		\ee
		Furthermore, if for probability measures $\nu_1$ and $\nu_2$ their Cauchy transforms are equal 
		$G_{\nu_1} = G_{\nu_2}$, then $\nu_1 = \nu_2$. 	
	\end{proposition}
	Notice, that the Stieltjes formula implies that the support of $\nu$ contains the set
	$$N =\{x\in\R:\, \Im G_\nu(x)\neq 0\}.$$
	
	The {moment transform} of $\X$ (or of its $\star$-distribution $\nu$) is 
	$$ M_\X(z) = M_\nu(z) = \int_\R (1-zt)^{-1}\, d\nu(t) $$
	for $z \in \C \setminus \R$. 
	It is clear that 
	\bel{for:G}
	G_\X(z)  = \frac{1}{z} \, M_\X\left(\frac{1}{z}\right).
	\ee
	\begin{proposition}[\cite{Vo93}]
		Let $\mu$ be a probabilistic measure on $[0,\infty)$ such that $\mu\left(\{0\}\right) \neq 1$.
		Then $M_\mu$ is injective on $i\,\C^+$ and the image $M_\mu(i\,\C^+)$ is an open set contained in~a~disc centered in~$(0,0)$ and of~diameter $1-\mu\left(\{0\}\right)$.
		Furthermore
		$$ M_\mu\left(i\C^+\right)\cap\R = \left( -\frac{1-\mu\left(\{0\}\right)}{2},\; \frac{1-\mu\left(\{0\}\right)}{2}\right).$$
	\end{proposition}
	\noindent So $\xi_\mu:\, M_\mu(i\C^+)\to i\C^+$, the inverse of $M_\mu$ is well defined.
	For a positive random variable $\X$ and its {$\star$-distribution} $\mu$ the { S--transform} is 
	$$ S_\X(z) = \frac{1+z}{z}\,\xi_\mu(z).$$
	\begin{proposition}[Theorem 4.5.3.23 from \cite{MS2017}]\label{T:S}
		If $\X$ and $\Y$ are free non-commutative positive random variables, then for $z$ in a neighbourhood of $0$:
		$$S_\X(z)\cdot S_\Y(z) = S_{\X\Y}(z).$$
	\end{proposition}
\subsection{Asymptotic freeness}
	Consider a family of random variables $\left( \X_1^{(n)},\ldots, \X_k^{(n)}\right)$ on a probabilisty space $(\A_n, \fif_n)$, $n=1,2,\ldots$. 
	We say that $\star$-distribution $\left(\mu_{\left( \X_1^{(n)},\ldots, \X_k^{(n)}\right)}\right)_{n\geq 1}$ converges to~$\mu$ when $n$ tends to infinity, 
	if for any polynomial in non-commutative variables ${P\in\C\left\langle x_1,\ldots, x_k\right\rangle}$:
	$$\int P(x_1,\ldots, x_k)\, d\mu_{\left( \X_1^{(n)},\ldots, \X_k^{(n)}\right)}(\x) 
	\stackrel{n\to\infty}{\longrightarrow} \int P(x_1,\ldots, x_k)\, d\mu(\x).$$
	If $\mu$ is a $\star$-distribution of $\left(\X_1,\ldots, \X_k\right)$, we say that $\left( \X_1^{(n)},\ldots, \X_k^{(n)}\right)$ {converges in distribution} to~$\left(\X_1,\ldots, \X_k\right)$. 
	Then
	$$\fif_n\left(P\left(\X^{(n)}_1,\ldots, \X^{(n)}_k\right)\right) 
	\stackrel{n\to\infty}{\longrightarrow} \fif\left(P(\X_1,\ldots, \X_k)\right)$$
	for any ${P\in\C\left\langle x_1,\ldots, x_k\right\rangle}$.
	If additionally $\X_1,\X_2,\ldots, \X_k$ are free, then we say that $\left(\X_1^{(n)}, \X_2^{(n)}, \ldots, \X_k^{(n)}\right)_{n\geq 1}$ is {asymptotically free}.

	{Empirical spectral distribution} of $N\times N$ random matrix $\x_N$ is a random measure
	$$P_N = \frac{1}{N}\left( \delta_{\lambda_1}+\delta_{\lambda_2} +\ldots + \delta_{\lambda_N}\right),$$
	where $\lambda_1, \lambda_2,\ldots, \lambda_N$ are eigenvalues of $\x_N$ (possibly multiple in algebraic meaning).

	Classical and powerful result on connection between random matrices and free random variables is due to~Voiculescu, \cite{Vo91}. Here we cite more general formulations from Chapter 4 in \cite{HP}.
	\begin{proposition}\label{voi}
		Consider $N\times N$ random matrices $X_N$ and $Y_N$ such that: 
		both $X_N$ and $Y_N$ have almost surely an asymptotic spectral distribution when $N\to\infty$; 
		$X_N$ and $Y_N$ are independent for $N\geq 1$; 
		$Y_N$ is unitarily invariant. 
		Then $X_N$ and $Y_N$ (as random variables in $(\A_N, 
		\fif_N)$, where $\A_N$ is the algebra of $N\times N$ matrices and $\fif_N = \tfrac{1}{N}\tr$) 
		are almost surely asymptotically free.
	\end{proposition}
	
	In \cite{Szp17} Szpojankowski used Theorem \ref{voi} to prove the Matsumoto-Yor property in a non--commutative setting. 
	Here we will adapt this approach to some extent in order to prove the non--commutative HV property. 
	The random matrix version of HV property is the starting point for this approach. It has been proved recently in~\cite{KP18} (Theorem 2.2.) and we cite it below.
	
	Let $M_N$ be the set of real symmetric $N\times N$ matrices. 
	By $M_N^+$ we denote the cone of positive definite symmetric matrices. 
	We say that random matrix $\x$ has the \textit{matrix--Kummer} distribution with parameters $a>\frac{N-1}{2}$, $b\in\R$ and $\Sigma \in M_N^+$, if it has the density:
	 \bel{eq:mk}
	 \mk_N\left(a,\, b,\, \Sigma\right) (dx) = \frac{1}{\Gamma _N(a)\,\Psi( a,\frac{N+1}{2}+b;\Sigma)}\,(\det x)^{a-\frac{N+1}{2}}\,(\det (\e+x))^{-(a+b)}\, e ^{-\langle\Sigma,\,x\rangle}\, I_{M_N^+}(x) \,dx,
	 \ee
	where $\Psi$ is the confluent hypergeometric function of the second kind of a matrix argument (see for~instance formula (2) in \cite{JJ85}), $\Gamma_N(x) = \pi^{N(N-1)/4}\prod_{j=1}^N (x+ (1-j)/2)$ and $\e$ is identity in $M_N^+$. Symbol ``$\<A,\,B\>$'' denots scalar product of matrices $A$ and $B$, so $\<A,\,B\> = \tr(AB)$.
	
	The Wishart distribution with parameters $b>\tfrac{N-1}{2}$ and $\Sigma \in M_N^+$, $\W(b,\Sigma)$, 
	has the density
	$$\W_N\left(b,\, \Sigma\right) (dx) = \frac{(\det \Sigma)^b}{\Gamma_r(b)}\, (\det y)^{b-\tfrac{N+1}{2}} \, e ^{-\langle\Sigma,\,x\rangle}\, I_{M_N^+}(x)\, dx.$$
	\begin{proposition}[Thm. 2.2. in \cite{KP18}]
	\label{tw1}
		 Let $X$ and $Y$ be two independent random matrices valued in $M_N^+$. 
		 Assume that $X$ has the matrix--Kummer distribution $\mk(a,\, b,\, c\,\e)$ and $Y$ has the Wishart distribution ${\W(a+b,\, c\,\e)}$, where $a>\frac{N-1}{2}$, $b>\frac{N-1}{2}-a$, $c>0$. 
		Then random matrices 
		$$U:=(\e+X)^{-1/2}\, Y\,(\e+X)^{-1/2},\qquad V:=(\e+U)^{1/2}\,X\,(\e+U)^{1/2}$$
		are independent. 
		Furthermore, $U\sim \mk(a+b,\, -b,\, c\e)$ and $V\sim \W(a,\,c\e)$.
	\end{proposition}
	\begin{remark}\label{rem:wishart}
	It is known that if we take a sequence of real Wishart matrices $(\y_n)$ with parameters $\alpha_n$ and $\Sigma_n = \lambda_n\e$,
	where $\lambda_n>0$, $\tfrac{2\lambda_n}{n}\rightarrow \lambda>0$ and $\tfrac{2\alpha_n}{n}\rightarrow \alpha$ 
	then the free Poisson distribution $\nu(1/\lambda, \alpha)$ is an almost sure weak limit of the empirical spectral distributions of $(\y_n)$.
	\end{remark}
	
	Probability measure	 $\nu(\lambda, \gamma)$ defined for $\lambda>0$ and $\gamma\geq 0$ by 
	$$\nu =\max\{0, 1-\lambda\}\,\delta_0 + \lambda\tilde\nu,$$
	where the measure $\tilde\nu$ is supported on the interval $\left(\gamma\, (1-\sqrt{\lambda})^2,\; \gamma\, (1+\sqrt{\lambda})^2\right)$ and has density
	$$\tilde\nu(dx) = \frac{1}{2\pi\gamma\,x}\, \sqrt{4\lambda\gamma^2 - \left(x- \gamma\,(1+\lambda)\right)^2}\, dx $$
	is called the free Poisson or Marchenko--Pastur distribution. 
	The parameter $\lambda$ is called the rate and $\gamma$ -- the jump size. 
	If $\X$ has free--Poisson $\star$-distribution with parameters $\gamma$ and $\lambda$, we denote it $\X\sim\mathit{f}\mathrm{Pois}(\gamma,\, \lambda)$.
	
%
\section{Asymptotic eigenvalue distribution of matrix--Kummer random matrix}\label{sec3}
In order to state the free HV property we need to find free counterpart of the classical Kummer distribution. 
In fact we will seek the asymptotic eigenvalue distribution of matrix--Kummer matrices. 

We recall a standard result on an eigenvalue distribution, which can be found in \cite{HP} (Prop. 4.1.3).
	\begin{proposition}\label{stwrozkladww}
		Let $\mathbf{Z}$ be symmetric $n\times n$ random matrix and let $g$ be its density with respect to the Lebesgue measure on $\R^{\tfrac{n\,(n+1)}{2}}$. 
		Assume that there exists a function $h:\R^n\to\R$ such, that 
		$$h\left(\lambda_1,\ldots, \lambda_n\right) = g\left(x_{ij},\; 1\leq i\leq j\leq n\right),$$
		 where $\lambda_1<\lambda_2< \ldots< \lambda_n$ are eigenvalues of symmetric $\x = (x_{ij})_{i,j=1,\ldots, n}$. 
		 \\
		Then the density of the vector of eigenvalues of~$\mathbf{Z}$ has the following form
		$$
		\frac{\pi^{n(n+1)/4}}{\prod_{j=1}^n\Gamma(j/2)}\,
		 h(\lambda_1, \ldots, \lambda_n)\,\prod_{i<j}|\lambda_i-\lambda_j| \,I(0<\lambda_1<\ldots<\lambda_n)
		.$$
	\end{proposition}
	
	\begin{definition}
		If $g$  is a density on $\R^n$ with respect to Lebesgue measure and 
		$$g(x) = C\cdot\prod_{1\leq i<j\leq n}|x_i-x_j|\,\prod_{i=1}^n\exp\left\{-\frac{n}{2}\,V_n(x_i)\right\}dx, $$
		where $x=(x_1, \ldots,x_n)$, then function $V_n:\R\to\R$ is called potential of $g$.
	\end{definition}
	
	Let ${\x_n\sim \mk_n(a_n, \,b_n,\, \gamma_n \e)}$, where $a_n>(n-1)/{2}$, $b_n\in\R$, $\gamma_n>0$.
	Then    by \eqref{eq:mk} the function $h$ defined in Prop.~\ref{stwrozkladww} has the form
	$$h_{(a_n,\,b_n,\,\gamma_n\e)} (\boldsymbol{\lambda}) = C\cdot \left(\prod_{i=1}^n \lambda_i\right)^{a_n-\frac{n+1}{2}}\left(\prod_{i=1}^n ( 1+\lambda_i) \right)^{-(a_n+b_n)} \exp\left\{-\gamma_n \sum_{i=1}^n \lambda_i \right\}\, I_{(0,\infty)^n}(\boldsymbol{\lambda}) 
	,
	$$
	where $\lambda_1<\lambda_2<\ldots<\lambda_n$ 
	and $\boldsymbol{\lambda} = (\lambda_1, \ldots, \lambda_n)$.
	\noindent
	From Prop.~\ref{stwrozkladww} it follows that eigenvalues of~$\x_n$ have joint density: 
	\bel{wwN}
		D\cdot \prod_{1\leq i<j\leq n} (\lambda_j - \lambda_i)\, e^{-\frac{n}{2} 
		\sum_{i=1}^n\left(\frac{ 2\gamma_n}{n}\lambda_i + \frac{2\beta_n}{n}\log(1+\lambda_i) 
		- \frac{2\alpha_n}{n}\log (\lambda_i)\right)}\,I(0<\lambda_1<\ldots<\lambda_n),
	\ee
	where $D$ is the normalizing constant, $\alpha_n = a_n -\tfrac{n+1}{2}>-1$ and $\beta_n = a_n + b_n\in\R$.

	From now on let $\beta_n>0$, $\alpha_n>-1$ and $\gamma_n>0$. 
	Moreover, let $V_n$ be the potential related to the density \eqref{wwN}, i.e.:
	\bel{eq:pot}
	V_n(x) = \frac{2\gamma_n}{n} x + \frac{2\beta_n}{n} \log(1+x) - \frac{2\alpha_n}{n}\log (x), \qquad x\in \R.
	\ee	 
	We also assume that ${2\gamma_n}/{n}\to \gamma>0$, ${2\beta_n}/{n}\to\beta\in \R$, $2a_n/{n}\to\alpha>0$, so that ${2\alpha_n}/{n} \to \alpha - 1$. 
	
	In the complex matrix case, one can set $\alpha_n = n\alpha$, $\beta_n = n\beta$ and $\gamma_n = n\gamma$ and then potential $V^{\mathrm{complex}}_n(x) = V^{\mathrm{complex}}(x)$ does not depend on $n$ and one can use classical results to get the~limit of empirical measure. 
	Here we use the following result from \cite{Fe08}: 
	
	If
	\begin{itemize}
		\item[(a)] for any $n\geq 1$: $V_n$ is  continuous;
		\item[(b)] there exist $a>0$ and $T>0$ such that
			$$ V_n(t)\geq (1+a)\,\log (1+t^2), \qquad \forall t\geq T;$$
		\item[(c)] there exists $V:\R^+\to \R$ such that $V_n$ uniformly converges to $V$ as $n\to\infty$ on compact subsets of $\R^+$,
	\end{itemize}
	then the random measure ${\hat{P}_n=\frac{1}{n}\sum_{i=1}^n \delta_{\lambda_i}}$ almost surely converges weakly to a probability measure $\mu_V$ as $n\to\infty$. 
	The measure $\mu_V$ minimizes the functional
	$$E_V(\mu) := \int\int_{\R^2} \log |s-t|^{-1}\,d\mu(s)\,d\mu(t) +\int V(s)\,d\mu(s)$$
	called the \textit{energy of a field with external potential $V$}.
	Properties of the measure $\mu_V$, called \textit{the equilibrium measure}, have been deeply analysed, see \cite{ST97} for instance. 
	The facts we use below also can be found in \cite{ST97} (see Theorems IV.1.11  and IV.3.1 there).
	
	If $V$ is convex on some closed interval $[A, B]$ or $xV'(x)$ is increasing on $(A,B)\subset[0,\infty]$
	then the support of the equilibrium measure $\mu_V$ is a closed interval $[a,b]\subset [A, B]$, where $a$ and $b$ are such that
	\bel{eq1}\left\{
		\begin{array}{ccl}
		\frac{1}{\pi}\int_a^b\frac{V'(x)}{\sqrt{(b-x)(x-a)}}\, dx &=& 0,\\
		\frac{1}{\pi}\int_a^b\frac{xV'(x)}{\sqrt{(b-x)(x-a)}}\,  dx &=& 2.
		\end{array}\right.
	\ee
	Then
	$$
		\mu_V(x) = \frac{1}{2\pi}\sqrt{(x-a)(b-x)}\, \,\, \,
		\mathrm{PV}\left(\frac{1}{\pi}\int_a^b \frac{V'(t)}{\sqrt{(t-a)(b-t)}} \frac{dt}{t-x}\right), \qquad x\in[a,b].
	$$
	\noindent
	For $V_n$ related to Kummer eigenvalues (see \eqref{eq:pot}) we have 
	\bel{pot}V(x) =\lim_{n\to\infty} V_n(x) = \gamma x + \beta \log (1+x) -(\alpha-1)\log(x).\ee
	That is conditions (a), (b), (c) are satisfied.
	Note that if $\beta \leq 0$ then $V$ is convex on $\R$ and if $\beta >0$ then $xV'(x)$ is increasing on $\R$. 
	So in both cases $[A,B] = [0,\infty)$ and \eqref{eq1} holds. 
	The system of equations~\eqref{eq1} can equivalently be written as:
	\bel{u2}
		\left\{
		\begin{array}{ccl}
		\gamma +\frac{\beta}{\sqrt{(a+1)(b+1)}} - \frac{\alpha-1}{\sqrt{ab}} &=& 0,\\
		\gamma\frac{a+b}{2} - \alpha+1 +\beta - \frac{\beta}{\sqrt{(a+1)(b+1)}} &=& 2,
		\end{array}\right.
	\ee
	 $a,b\in[0,\infty)$. 
	To see this, one can use Cauchy's residual theorems to calculate the integrals, which appear in \eqref{eq1}.
	
	Moreover, in this case the equilibrum measure is 
	\bel{muV}
		\mu_V(x) = \frac{1}{2\pi}\sqrt{(x-a)(b-x)}\left( \frac{\alpha-1}{x\sqrt{ab}} 
		- \frac{\beta}{(1+x)\sqrt{(a+1)(b+1)}}\right) I_{(a,b)}(x).
	\ee
	Note that the measure $\mu_V$ depends on $\gamma$ through \eqref{u2}.
	 The distribution defined in~\eqref{muV} will be called the \textit{free--Kummer} distribution with parameters $\alpha$, $\beta$, $\gamma$.
	 We will write ${\X\sim\mathit{f}\K(\alpha, \beta, \gamma)}$ if $\mu_V$ is $\star$-distribution of $\X$.

We sum up preceding calculations in the remark below.	
	\begin{remark}\label{rem:kummer}
	If $\x_n\sim\mk(a_n,\, b_n,\,c_n)$, where $\frac{n-1}{2}<a_n$ and $\tfrac{a_n}{n}\to a/2>0$, $\tfrac{b_n}{n}\to b/2\in\R$ and $\tfrac{c_n}{n}\to c/2>0$, then the limiting spectral distribution of $(\x_n)_n$ is free--Kummer $\mathit{f}\K(a,\,a+b,\,c)$.
	\end{remark}
	%
	
	\begin{lemma}\label{cauchykummer}
		The Cauchy transform of free--Kummer $\mathit{f}\K(\alpha, \beta, \gamma)$  is 
		\bel{cauchy}
			G(z)= \frac{1}{2}\left(\gamma - \frac{\alpha-1}{z} + \frac{\beta}{1+z} 
			+\sqrt{(z-a)(z-b)}\,\left[ \frac{\beta}{(1+z)\sqrt{(a+1)(b+1)}}-\frac{\alpha-1}{z\sqrt{ab}} \right] \right)		
		\ee
	\end{lemma}

	Our approach to proving Theorem~\ref{wprost} (see  Lemma~\ref{limits}) requires the largest eigenavlue of Kummer matrix to be asymptotically a.s. bounded. 
	For this reason we introduce the large deviation principle (LDP).
	We say that the LDP with a rate function $I$ and speed $n$ holds for a sequence of measures $(\nu_n)$, 
	if for any Borel set $\Gamma$
	$$
		-\inf\{I(x) : x\in \Gamma^\circ\}\leq \liminf_{n\to\infty} n\log \nu_n(\Gamma)\leq\limsup_{n\to\infty}
		 n\log \nu_n(\Gamma)\leq -\inf\{I(x) : x\in \bar\Gamma\},
	$$
	\noindent where $\Gamma^\circ$ is the interior and $\bar\Gamma$ is the closure of $\Gamma$. 
	Function $I$ is a {\bf good rate function}, if for~any $\alpha\in\R$ the set $\{x: I(x)\leq \alpha\}$ is compact.	
	\begin{proposition}\label{Largest}
		The largest eigenvalue $\lambda_{\mathrm{max}}$ of the Kummer matrix $\mk(\alpha_n,\beta_n, \gamma_n)$  converges almost 
		surely to $b$ (right end of the support) and satisfies the LDP on $\R^{+*}$ with speed $n$ and the good rate function
		$$ I^*_{\alpha, \beta, \gamma}(t) = 
		\left\{
		\begin{array}{ll}\displaystyle 
			\int_b^t \tfrac{1}{2}\sqrt{(x-a)(x-b)}\,\left(\tfrac{\alpha -1 }{x\sqrt{ab}} - \tfrac{\beta}{(1+x)\sqrt{(a+1)(b+1)}}\right)dx, & t>b\\\displaystyle 
			+\infty, & \mathrm{otherwise}
		\end{array}\right.$$
	\end{proposition}

	\begin{corollary}\label{ldp_wniosek}
		The largest eigenvalue of matrix Kummer random matrix is asymptotically almost surely  bounded.
	\end{corollary}

\section{Freeness property of free--Poisson and free--Kummer free variables}\label{sec4}
Now, we are ready to state HV property for non--commutative random variables. 
\begin{theorem}\label{wprost}
	Let $(\mathcal A, \fif)$ be a $C^\star$-probability space. 
	Assume that $\X,\Y\in\A$ be free random variables and
	${\X\sim \mathit{f}\K (\alpha,\, \alpha+\beta,\, \gamma)}$, $\Y\sim\mathit{f}\mathrm{Pois} (1/\gamma,\, \alpha+\beta)$, with $\alpha,\gamma>0$ and $\beta>-\alpha$. 
	Let
	$ \U : = (\i+\X)^{-1/2}\,\Y\,(\i+\X)^{-1/2}$ and $\V :=(\i+\U)^{1/2}\,\X\,(\i+\U)^{1/2}$.
	Then $\U$ and $\V$ 
	are free. Moreover, $\U\sim\mathit{f}\K(\alpha+\beta,\,\alpha,\,\gamma)$ and $\V\sim\mathit{f}\mathrm{Pois}(1/\gamma,\,\alpha)$. 
\end{theorem}
 
To prove Theorem~\ref{wprost}, we need a technical result related to convergence in probability, which is a generalization of a lemma from \cite{Szp17}.

\begin{lemma}\label{limits}
	Let $(\u_N)_{N\geq 1},\; (\z_N)_{N\geq 1}$ be two independent sequences of random matrices on~probabilistic space $(\Omega,\,\mathcal{F},\, \P)$, where $\u_N$, $\z_N$ are $N\times N$ matrices for every $N$.
	Suppose that $\u_N$ and $\z_N$ have almost surely weak limits of their sequences of the empirical spectral distributions, 
	$\mu$ and $\nu$, respectively. 
	Also suppose that the smallest eigenvalue of $\z_N$ is asymptotically almost surely larger than a constant $A>0$ and 
	that the largest eigenvalue of $\z_N$ is asymptotically almost surely smallert than a constant $B>0$. 
	Let $(\mathcal A,\, \fif)$ be $C^\star$-probability space.
	
	Assume that there exist $\U,\, \Z\in\mathcal A$ such that $\U$ and $\Z$ are free and $\U\sim \mu$, $\Z\sim \nu$.	
	Then for any complex polynomial $Q\in \C\left\<x_1, x_2, x_3\right\>$ in three non-commuting variables and 
	for any $\epsilon>0$ we have
	$$\P\left( \left\vert\frac{1}{N}\tr\left[Q\left(\u_N,\, \z_N,\, \z_N^{-1}\right)\right] - \fif\left[Q\left(\U,\, \Z,\, \Z^{-1}\right)\right]\right\vert>\epsilon\right)\stackrel{N\to\infty}{\longrightarrow} 0. $$
\end{lemma}
In \cite{Szp17} this lemma was proved in special case of $\U$ with free GIG distribution, and $\Z$ with free Poisson distribution. 
Although our formulation is more general the proof remains the same, so we skip it.
\begin{corollary}\label{wn1}
	Let $(\mathcal A, \fif)$ be a $C^\star$-probability space.
	Assume that there exist $\X, \Y\in\mathcal A$ such that $\X$ and $\Y$ are free, 
	$\X\sim \mathit{f}\K(a,\, a+b,\,c)$ and $\Y\sim \mathit{f}\mathrm{Pois}(1/c,\, a+b)$.	
	Let $(\x_N)_{N\geq 1},\; (\y_N)_{N\geq 1}$ be two sequences of~random matrices, such that
	 $\x_N\sim \mk(a_N,\,b_N,\, c_N\e)$ and ${\y_N\sim\W(a_N+b_N,\, c_N\e)}$ are independent for each $N$ and $2\, a_N/N\to a$, $2\, b_N/N\to b$, $2\,c_N/N\to c$. 	
	Then for any complex polynomial $Q\in \C\left\<x_1, x_2, x_3\right\>$ in three non-commuting variables and 
	for any $\epsilon>0$ we have
	\bel{lim1}
		\P\left( \left\vert\frac{1}{N}\tr\left[Q\left(\e+\x_N,\, (\e+\x_N)^{-1},\, \y_N\right)\right] 
		- \fif\left[Q\left(\i+\X,\, (\i+\X)^{-1},\, \Y\right)\right]\right\vert>\epsilon\right)\stackrel{N\to\infty}{\longrightarrow} 0. 
	\ee
\end{corollary}

\section{Free regression characterization of HV type}\label{sec5}
\subsection{Conditional expectation in a non-commutative  probability space}
	In the next subsection we will formulate free version of the following characterization theorem from \cite{PW18}, which holds in the classical 
probability setting.
	\begin{proposition}\label{reg}
		Let $X$ and $Y$ be independent, positive, non-degenerate (commutative) random variables, 
		such that ${\E X <\infty}$, $\E Y <\infty$ and $\E X^{-1} <\infty$. 
		Let $U:={Y}/({1+X})$, $V:= X\left(1+U\right)$ and  assume that there exist real constants $\alpha$ and $\beta$ such that
		$\E\left(V|U\right) = \alpha$ and 
		$ \E\left(V^{-1}|U\right) = \beta. $
		Then $\alpha\, \beta>1$ and there exists a constant $c>0$ such that
		$$
		 X\sim\mathcal{K}\left( \frac{\alpha\beta}{\alpha\beta -1},\,
		 c - \frac{\alpha\beta}{\alpha\beta -1},\,
		 \frac{\beta}{\alpha\beta -1}\right) \qquad\mbox{and}\qquad  Y\sim \mathcal{G}\left(c,\,\frac{\beta}{\alpha\beta - 1}\right).
		$$
	\end{proposition}
	We will recall the definition of non-commutative conditional expectation following \cite{BB04, Take79}. 
	Let $(\mathcal A, \fif)$ be a $W^\star$--probability space.
	Let $\mathcal B\subset \mathcal A$ be a von Neumann subalgebra of $\mathcal A$. 
	Then there exists a unique faithful, normal projection $\fif(\cdot|\mathcal B): \A\rightarrow\mathcal B$ such that 
	$\fif(\fif(\cdot|\mathcal B))=\fif(\cdot)$. 
	We call it a non-commutative conditional expectation from $\mathcal A$ to $\mathcal B$ with respect to $\fif$ (see \cite{Take79}, Vol I p. 332). 
	The conditional expectation of a self--adjoint element $\X\in\mathcal A$ is a unique self--adjoint element of $\mathcal B$. 
	
We cite, following \cite{BB04}, two 
important properties of  a non-commutative conditional expectation $\fif(\cdot\vert\mathcal B)$, which will be used in the proof of Theorem \ref{free_char}.	
	\begin{enumerate}[label=\Roman{*}.]
	\item\label{p1} If random variables $\U,\, \V\in\mathcal A$ are free, then $\fif\left(\U\vert \V\right) = \fif\left(\U\right) \i$;
	\item \label{p2} If $\X\in \mathcal A$, $\Y\in\mathcal B$, then $\fif\left(\X\Y\right) = \fif\left(\fif\left(\X\vert\mathcal B\right)\Y\right).$
	
	\end{enumerate}

\subsection{The characterization theorem}

\begin{theorem}\label{free_char}
	Let $(\A, \phi)$ be a non--commutative $W^\star$--probability space. 
	Let $\X\in\A$ and $\Y\in\A$ be self--adjoint, positive, free, compactly supported and non-degenerate random variables. 
	Define
		$$\U := (\i+\X)^{-1/2}\, \Y\, (\i+\X)^{-1/2}\qquad\mathit{and}\qquad \V := (\i+\U)^{1/2}\,\X\,(\i+\U)^{1/2}$$ 
	and assume that there exist constants $\bar\alpha,\bar\beta>0$ such that
	\bel{jed} \fif\left(\V\,|\,\U\right) = \bar\alpha\, \i,\ee
	\bel{dwa}\fif\left(\V^{-1}\,|\,\U\right)  = \bar\beta\, \i.\ee	 
	 Then $\bar{\alpha}\bar{\beta}>1$ and there exists $a>0$ such that $\X\sim \mathit{f}\K(\bar\alpha\gamma,\, a\gamma,\,\gamma$),
	 $\Y\sim\mathit{f}\mathrm{Pois}(1/\gamma,\,a\gamma)$, 
	where ${\gamma={\bar\beta}/\left({\bar\alpha\bar\beta-1}\right)}$.
\end{theorem}

\section{Concluding remarks}
We have proved that free Poisson and free Kummer are the only probability distributions which maintain freeness of non-cummutative random variables when transformed by the following mapping 
\small
$$
	(x,\,y)\mapsto \left( (\i + x)^{-1/2}\, y\, (\i+x)^{-1/2},\;
	\left[ (\i + x)^{-1/2}\, y\, (\i+x)^{-1/2}\right]^{1/2}\, 
	x  \left[ (\i + x)^{-1/2}\, y\, (\i+x)^{-1/2}\right]^{1/2} \right). $$
\normalsize
As has been mentioned in \cite{KP18}, in the matrix setting one can consider a different transformation:
\small
$$
	(x,\,y)\mapsto \left((\e+x+y)^{1/2}\,(\e+x)^{-1}\,(\e+x+y)^{1/2} -\e,\; x+y - \left[(\e+x+y)^{1/2}\,(\e+x)^{-1}\,(\e+x+y)^{1/2} -\e\right]\right) $$
\normalsize	
that preserves independence for Kummer and Wishart random matrices. It is still unknown if a related characterization holds for random matrices as well as if its free counterpart is true. 

An open question of a broader nature would be: does every independence characterization of random matrices has its analogon in free probability. Examples that have been studied suggest that answer could be positive. If so, then how to find it?

\section{Proofs of Sec. \ref{sec3}}
\subsection{Proof of Lemma \ref{cauchykummer}}
	\begin{proof}
	It is enough to show that
	\begin{enumerate}
		\item 
		$$ \int_a^b\sqrt{(x-a)(b-x)} \,\frac{1}{z-x\,}\frac{1}{x}\,dx = \frac{\pi}{z}\left(z-\sqrt{ab} - \sqrt{(z-a)(z-b)}\right),$$
		\item 
		$$ \int_a^b\sqrt{(x-a)(b-x)} \,\frac{1}{z-x}\,\frac{1}{1+x}\,dx = \frac{\pi}{1+z}\,\left((z+1) - \sqrt{(a+1)(b+1)} - \sqrt{(z-a)(z-b)}\right).$$
	\end{enumerate}
	These equalities can be obtained using the Cauchy Residual Theorem.
	The square root above denotes its main branch. 
	So 	$\sqrt{\C^-}\subset \C^-$ and~${\sqrt{\C^+\cup \R} \subset \C^+\cup\R^+\cup \{0\}}$.
	Another choice of branch is not possible, since $\mu$ has a compact support and $\displaystyle \lim_{x\to\infty,\, \eps\to 0} \Im G_\mu(x + i\eps) $ has to be equal to $0$ due to Thm.~\ref{Stieltjes} (Stieltjes' Inversion Formula).
\end{proof}

\subsection{Proof of Proposition \ref{Largest}}
\begin{proof}
	We repeat the reasoning from the proof of Theorem 2 in \cite{Fe06}. 
	
	Let $$g(x) := \int_a^b \log|x-t|^{-1}\,\mu(dt) + \tfrac{1}{2}V(x) + \tfrac{1}{2}\int_a^bV(t)\, \mu(dt),$$
	\noindent 
	where $V$ is defined in \eqref{pot}.
	Notice that $I^*_{\alpha, \beta, \gamma}(t) = \int_b^t g'(x) \,dx$ for all $t>b$. 
	Since \vspace*{0.3cm}
	$$g'(x) = \tfrac{1}{2}V'(x) + \int_a^b \tfrac{1}{t-x} \,d\mu(t) = \tfrac{1}{2}\sqrt{(x-a)(x-b)}\left(\tfrac{\alpha -1 }{x\sqrt{ab}} - \tfrac{\beta}{(1+x)\sqrt{(a+1)(b+1)}}\right)$$ 
	\noindent is positive ($\gamma>0$), then $g$ is increasing on $(b, \infty)$. 
	The remaining part of the proof follows from the arguments in~Section 4.2 in \cite{Fe06}.
\end{proof}

\section{Proofs of Sec. \ref{sec4}}
\subsection{Proof of Corollary \ref{wn1}}
\begin{proof}
	The sequences of empirical spectral distribution of matrices $(\x_N)_N$ and $(\y_N)_N$ almost surely have their weak limits: 
	$\mathit{f}\K(a,\,a+b,\,c)$ and $\mathit{f}\mathrm{Pois}(1/c,\,a+b)$. 
	Eigenvalues of $\e+\x_N$ are greater than $1$.
	It implies that the support of the weak limit of empirical spectral measure of sequence $(\e+\x_N)_{N=1,2,\ldots}$ is
	separated from $0$. 
	Also the largest eigenvalue is asymptotically almost surely bounded as $N\to\infty$ (Cor.~\ref{ldp_wniosek}). 
	Then the result follows from Lemma~\ref{limits}.
\end{proof}

\subsection{Proof of Thm. \ref{wprost}}
\begin{proof}
	We want to show that the algebras generated by $\U$ and $\V$ are free.
	
	Let us take a sequence $(\y_n)_{n\geq 1}$ of $n\times n$ Wishart matrices with parameters $\alpha_n+\beta_n$ and $\gamma_n \e$, such that $\alpha_n>\tfrac{n-1}{2}$, $\alpha_n + \beta_n>\tfrac{n-1}{2}$, $\gamma_n>0$ and 
$2\alpha_n/n\to \alpha>0$, $2\beta_n/n\to \beta$ and $2\gamma_n/n\to \gamma>0$. 
	Moreover let $(\x_n)_{n\geq 1}$ be  a sequence of Kummer matrices with parameters $\alpha_n$, $\beta_n$ and $\gamma_n$. 
	Assume that $(\x_n)$ and $(\y_n)$ are independent. 
	Then since matrix--Kummer and Wishart matrices are unitarily invariant and both have almost sure limiting eigenvalue distributions, 
	they are asymptotically free (Thm. \ref{voi}). 
	It means that if $\mathcal A_n$ is the algebra of random matrices of size $n\times n$ with integrable entries 
	with the state $\fif_n(a) = 1/n\tr (a)$ on $\A_n$, 
	then for any polynomial in~two non--commuting variables $P\in\C\left\<x_1,\, x_2\right\>$ we have almost surely
	$$\lim_{n\to\infty} \fif_n(P(\x_n,\, \y_n)) = \fif (P(\X,\, \Y)),$$
	where $\X$ and $\Y$ are as in the statement of the theorem. 
	
	By the HV property for random matrices it follows that for any $n\in\N_+$
	$$\u_n := (\e+\x_n)^{-1/2}\, \y_n\,(\e+\x_n)^{-1/2}\qquad\mathrm{and}\qquad\v_n := (\e+\u_n)^{1/2}\, \x_n\, (\e+\u_n)^{1/2}$$ 
	are independent, $\u_n\sim\mk(a_n+b_n,\, -b_n, \, \gamma_n\e)$ and $\v_n\sim \W(a_n,\, \gamma_n\e)$. 
	They are also almost surely asymptotically free. 
	Let $\tilde\U, \tilde\V$ be the limiting pair of non--commuting free random variables. 
	
	Then for any polynomial $P\in\C\left\<x_1,\, x_2\right\>$, there exists $Q\in\C\<x_1, x_2, x_3\>$ such, that we have 
	\bel{w1}
	\lim_{n\to\infty}\fif_n\left(Q\left(\e+\x_n,\, (\e+\x_n)^{-1},\, \y_n\right)\right)=\lim_{n\to\infty} \fif_n(P(\u_n,\, \v_n)) =\fif(P(\tilde\U,\, \tilde\V)).\ee
	
	From Corollary \ref{wn1} we know that $\fif_n\left(Q\left(\e+\x_n,\, (\e+\x_n)^{-1},\, \y_n\right)\right)$ converges in probability to
	$$\fif\left(Q\left(\i+\X,\, (\i+\X)^{-1},\, \Y\right)\right).$$
	However, by \eqref{w1}  the sequence $(\fif_n\left(Q\left(\e+\x_n, (\e+\x_n)^{-1}, \y_n\right)\right))_{n\geq1}$ has almost sure limit and thus we have
	$$\lim_{n\to\infty}\fif_n\left(Q\left(\e+\x_n,\, (\e+\x_n)^{-1},\, \y_n\right)\right)= \fif\left(Q\left(\i+\X,\, (\i+\X)^{-1},\, \Y\right)\right) = \fif\left(P\left(\U,\, \V\right)\right),$$
	where the last equality follows from the relation between $P$ and $Q$.
	
	Thus joint moments of $(\U, \V)$ and $(\tilde \U,\tilde \V)$ are the same. 
	Since $\tilde \U$ and $\tilde \V$ are free, then $\U$ and $\V$ are also free (recall that freeness is defined by joint moments).
	
	From Remark \ref{rem:wishart} and Remark \ref{rem:kummer} we have that $\U\sim\mathit{f}\K(\alpha+\beta,\,\alpha,\,\gamma)$ and $\V\sim\mathit{f}\mathrm{Pois}(1/\gamma, \, \alpha)$.
\end{proof}

\section{Proofs of Sec. \ref{sec5}}
\subsection{Proof of Theorem \ref{free_char}}
\begin{proof}
\textbf{Step 1.}
	First step is to show that $\Y$ is free--Poisson random variable.
	
	For any $k\geq 0$ we multiply both sides of \eqref{jed} by $\U^k$ and take expectation $\fif$. 
	Due to the property (II) we have
	$$\fif\left(\U^k(\i+\X) + \U^{k+1}(\i+\X) - \U^k - \U^{k+1}\right) =\bar \alpha\,\fif\left(\U^k\right).$$
	Similarly, if we multiply Eq. \eqref{dwa} by $(\i + \U)\U^k$ and take expectation, we have
	$$
	\fif\left(\U^k\X^{-1}\right) = \bar\beta\fif\left(\U^k (\i+\U)\right).
	$$
	We obtain a system of recursive equations holding for any $k\geq 1$
	\bel{rek}
		\begin{array}{ccl}
		\beta_{k-1}+\beta_k &=& (\bar\alpha + 1 )\, \alpha_{k} + \alpha_{k+1},\\
		\gamma_k &=& \bar\beta\,(\alpha_k + \alpha_{k+1}),
		\end{array}
	\ee
	where for $k\geq 0$: 
	\begin{equation*}
	\begin{split}
		\alpha_k &= \fif\left(\left[(\i+\X)^{-1}\,\Y\right]^k\right) =  \fif\left(\left[\Y\,(\i+\X)^{-1}\right]^k\right),
		\\ \beta_k &= \fif\left(\Y\left[(\i+\X)^{-1}\,\Y\right]^k\right) =  \fif\left(\Y\left[\Y\,(\i+\X)^{-1}\right]^k\right),
		\\ \gamma_k &= \fif\left(\X^{-1}\left[(\i+\X)^{-1}\,\Y\right]^k\right) =  \fif\left(\X^{-1}\left[\Y\,(\i+\X)^{-1}\right]^k\right).
	\end{split}
	\end{equation*} 
	For instance, notice that 	
	$\alpha_0 = 1$ and $\beta_0 = \fif(\Y)$.
	For $z$ from neighbourhood of $0$ we can define 
	$${A(z): = \sum_{n=0}^\infty z^n \alpha_n},\qquad 
	B(z): = \sum_{n=0}^\infty z^n \beta_n,\qquad 
	C(z): = \sum_{n=0}^\infty z^n \gamma_n.
	$$ 
	From  Eq. \eqref{rek} we obtain
	\begin{enumerate}
		\item[(d)] 
		$$B(z) +\frac{B(z) - \beta_0}{z} = \frac{\bar\alpha +1}{z}\, (A(z) -1) + \frac{A(z) - z\alpha_1 - 1}{z^2},$$
		\item[(e)] 
		$$C(z) =\bar\beta\left(A(z) +\frac{A(z) - 1}{z}\right).$$
	\end{enumerate}
	Also, denote $r:=r_\Y$ a $r$-transform of $\Y$ and $$D(z) : = \sum_{n=0}^\infty z^n \delta_n,$$ where 
	$$\delta_n= \fif\left( (\i+\X)^{-1}\left[\Y\,(\i+\X)^{-1}\right]^n\right).$$
	Then these three relations hold:
	\begin{itemize}
		\item[(a)] $A(z) = 1+ zD(z)\,r\left(zD(z)\right)$: 
		\\From \eqref{wzorBLS} it follows that for any $n\geq 1$ we have
		$$
		 \alpha_n =\displaystyle \fif\left(\left[ (\i + \X)^{-1}\, \Y\right]^n\right) 
		 = \sum_{k=1}^n\ka_k \sum_{i_1+\ldots+i_k = n-k}\delta_{i_1}\ldots\delta_{i_k},
		$$ 
		where $\ka_k :=\ka_k(\Y)$ is $k$-th free cumulant of $\Y$.

		So
		\begin{equation*}\begin{split}
			 A(z) &= \sum_{n=1}^\infty  z^n \, \sum_{k=1}^n\ka_k \sum_{i_1+\ldots+i_k = n-k}\delta_{i_1}\ldots\delta_{i_k} +1 \\
			&=\sum_{k=1}^\infty \ka_k \, z^k \sum_{n=k}^\infty  z^{n-k}\sum_{i_1+\ldots+i_k = n-k}\delta_{i_1}\ldots\delta_{i_k} +1\\
			&= \sum_{k=1}^\infty \ka_k \, z^k D^k(z) + 1\\ 
			&= 1+zD(z)\,r\left(zD(z)\right).
		\end{split}
		\end{equation*}

	\item[(b)] $B(z) = A(z)\, r\left(zD(z)\right)$:\\
		From \eqref{wzorBLS} it follows that for any $n\geq 0$ we have
		$$\beta_n = \fif\left(\Y\left[ \Y(\i + \X)^{-1}\right]^n\right) 
		=\displaystyle \sum_{k=1}^{n+1}\ka_k 
		\sum_{i_1+\ldots+i_k = n-k}\alpha_{i_1}\delta_{i_2}\ldots\delta_{i_k}.$$
		 Then
		\begin{equation*}\begin{split}
			B(z) &=\sum_{n=1}^\infty z^n \,\sum_{k=1}^{n+1}\ka_k  \sum_{i_1+\ldots+i_k = n-k}\alpha_{i_1}\delta_{i_2}\ldots\delta_{i_k} \\
			&= \sum_{k=1}^\infty \ka_k \, z^k \sum_{n=k-1}^\infty z^{n-k}\sum_{i_1+\ldots+i_k = n-k}\alpha_{i_1}\delta_{i_2}\ldots\delta_{i_k} \\
			&= \sum_{k=1}^\infty \ka_k \, z^k\, A(z) D^{k-1}(z)\\ 
			&= zA(z)\,r\left(zD(z)\right).
			\end{split}
		\end{equation*}

		\item[(c)] $C(z)\left[z\,r(zD(z))-1\right] = A(z)-1-\gamma_0$:
		\\
		Note that for $n\geq 1$ we have
		\nobel\nonumber \gamma_n =\fif\left( \X^{-1} \left[(\i+\X)^{-1} \,\Y \right]^n\right) = \fif\left(\X^{-1} \Y \left[(\i+\X)^{-1}\,\Y\right]^{n-1}\right) - \fif\left(\left[(\i+\X)^{-1}\,\Y\right]^n\right). \ee
		Again Eq.~\eqref{wzorBLS} implies
		\nobel\nonumber\begin{array}{ccl}
		\displaystyle \fif\left(\X^{-1} \Y \left[(\i+\X)^{-1}\,\Y\right]^{n-1}\right) &=& \ka_1\gamma_{n-1} + \ka_2\, (\delta_0 \gamma_{n-2} + \delta_1\gamma_{n-3} +\ldots+\delta_{n-2}\gamma_0)\, +\\
		 &&+\ldots +\ka_n\, \delta_0^{n-1}\gamma_0 \\
		&=& \displaystyle\sum_{k=1}^n \ka_k\displaystyle\sum_{i_1+\ldots+i_k = n-k} \gamma_{i_1}\delta_{i_2}\cdot\ldots\cdot\delta_{i_k}.
		\end{array}
		\ee
		
		If we multiply $\gamma_n$ by $z^n$ and sum over $n=0,1,\ldots$, we have:
		\nobel\nonumber\begin{split}
		C(z) &= \gamma_0 + \sum_{n=1}^\infty z^n\sum_{k=1}^n \ka_k\sum_{i_1+\ldots+i_k = n-k} \gamma_{i_1}\delta_{i_2}\cdot\ldots\cdot\delta_{i_k} - \sum_{n=1}^\infty z^n \alpha_n\\
		&= \gamma_0 +\sum_{k=1}^\infty\ka_k\, z^k\sum_{n=k}^\infty z^{n-k}\sum_{i_1+\ldots+i_k = n-k} \gamma_{i_1}\delta_{i_2}\cdot\ldots\cdot\delta_{i_k} - (A(z) -1)\\
		&=\gamma_0 +  z C(z)\sum_{k=1}^\infty \ka_k\, z^{k-1} D(z)^{k-1} - A(z) +1\\
		&= \gamma_0 +1 + z\,C(z)\, r(zD(z)) - A(z).
		\end{split}
		\ee
	\end{itemize}
	
	To determine distribution of $\Y$ we will solve the system of equations (a)-(e) with respect to~$r$. 
	We rewrite the system here:
	\begin{itemize}
		\item[(a)] $A(z) = 1+ zD(z)\,r(zD(z))$,
		\item[(b)] $B(z) = A(z) \, r(zD(z))$,
		\item[(c)] $\displaystyle C(z)\left[z\,r(zD(z))-1\right] = A(z)-1-\gamma_0$,
		\item[(d)] $\displaystyle B(z) +\frac{B(z) - \beta_0}{z} = (\bar\alpha +1)\, \frac{A(z) -1}{z} +	 \frac{A(z) - z\alpha_1 - 1}{z^2}$,
		\item[(e)] $\displaystyle C(z) =\bar\beta\left(A(z) +\frac{A(z) - 1}{z}\right)$.
	\end{itemize}
	Firstly, we multiply Eq. (e) by $z\,r(zD(z)) - 1$ and then we plug it into Eq.~(c). 
	We obtain	
	$$
	\bar\beta^{-1}\,\left(  A(z) - 1-\gamma_0\right) = \left( A(z) 
	+ \frac{A(z) -1}{z}\right) \left(z\, 	r(zD(z)) - 1\right).
	$$
	Let $h(z) := zD(z) r(zD(z))$. 
	Then (a) can be written in terms of $h$ as: $A(z) = 1+h(z)$.
	From this and above equation we get:
	$$
	\bar\beta^{-1} \left(  h(z)-\gamma_0\right) = 
	\left( 1+h(z) + \frac{1}{z}h(z)\right) \left(\frac{h(z)}{D(z)} - 1\right).
	$$
	Then after simplification
	\bel{wst1}
	\left(h(z) +\frac{1}{z}h(z)\right)\frac{h(z)}{D(z)}=\bar\beta^{-1}\, h(z)-\bar\beta^{-1}\, 
	\gamma_0- \frac{h(z)}{D(z)} + 1+h(z)+\frac{h(z)}{z}.
	\ee
	On the other hand we can plug $B(z)$ from (b) into (d) and then multiply it by $z$. 
	Now we have
	\bel{odw1}\nonumber(1+z)A(z)\,r(zD(z)) -\beta_0 -\bar\alpha\,(A(z) - 1) + \alpha_1+1 
	= A(z) +\frac{A(z) -1}{z}.\ee 
	In terms of $h(z)=zD(z)\, r(zD(z)) = A(z) -1$ it reads as
	\bel{wst2}
		\left(h(z) +\frac{1}{z}h(z)\right)\frac{h(z)}{D(z)}= 
		h(z) + \frac{h(z)}{z} +\beta_0 +\bar\alpha h(z)-\alpha_1 -\frac{h(z)}{D(z)} - \frac{h(z)}{zD(z)}.
	\ee
	
	Comparing the left-hand sides of \eqref{wst1} and \eqref{wst2}, we arrive at
	$$
	\bar\beta^{-1}\,h(z)-\bar\beta^{-1}\,\gamma_0+1= \beta_0+\bar\alpha \,h(z)-\alpha_1-\frac{h(z)}{zD(z)}.
	$$
	Which can equivalently be written as
	$$
	h(z) =\left({\beta_0 - \alpha_1 - 1 +\bar\beta^{-1}\, \gamma_0}\right)\,\left(\bar\beta^{-1} -\bar \alpha +\frac{ 1}{zD(z)}\right)^{-1}.
	$$
	Since 
	the $r$-transform of $\Y$ is analytic near $0$ and $\lim_{z\to 0 } zD(z) = 0$ (see definition of $D$), then for $y$ close to $0$ we have:
	\begin{equation*}
		r(y) =\left({\beta_0 - \alpha_1 - 1 +\bar\beta^{-1}\, \gamma_0}\right)\, \left({y\left(\bar\beta^{-1} - \bar\alpha\right) + 1}\right)^{-1}=\frac{\beta_0}{  1-y\left(\bar \alpha-\bar\beta^{-1}\right)} .
	\end{equation*}
	The last equality holds due to the fact, that $\alpha_1+1 = \gamma_0/\bar\beta$.   
	Eventually, we have
	\bel{rtransy}r(y) = \frac{a}{1-by},
	\ee
	where $a = \beta_0$ and $b=\bar\alpha - 1/\bar\beta$.  
	It is an $r$-transform of free-Poisson distribution with parameters $({\bar\alpha\bar\beta-1})/{\bar\beta}$ and ${a\bar\beta}/({\bar\alpha\bar\beta-1})$.
	 
	\textbf{Step 2.} 
	To recover Cauchy transform of $\U$ we use (a), (b) and (d). From (a) and formula for $r$ we can deduce  that for $z$ in a neighbourhood of $0$
	$$zD(z) = \frac{A(z) -1}{(A(z) -1) b +a}.
	$$ 
	Plugging it into (b) we 
get $B(z) = A(z)( (A(z)-1)b + a)$. We plug that last in (d), which gives quadratic equation for $A$:
	$$
		A^2(z) z (1+z) (\bar\alpha\bar\beta-1) + A(z)\left[z(1+z)(\bar\beta a-\bar\beta\bar\alpha +1)
		 - \bar\beta(1+z+\bar\alpha z)\right]+ \bar\beta +z\tilde{c}=0,
	$$
	where $\tilde{c}=\bar\beta(1+\bar\alpha - a +\alpha_1)$.
	Solution has the following form:
	\footnotesize
	\begin{multline*}  
		A(z)= \frac{\bar\beta}{2 (-1+\bar\alpha \bar\beta)\, z\, (1+z)}\bigg(1+\left(-\bar\beta^{-1}+1-a +2 \bar\alpha \right)
		\,z+\left(-\bar\beta^{-1}-a+\bar\alpha \right)\, z^2+\bigg. \\
		\left.+\bar\beta^{-1}\sqrt{ -4\bar \beta (-1+\bar\alpha \bar\beta) 
		z (1+z) \left[1+(1-a+\bar\alpha+\alpha_1) z\right]+\left[z^2(1+\bar\beta(a-\bar\alpha)) 
		+ z(1+\bar\beta(a - 1 - 2\,\bar\alpha)) - \bar\beta\right]^2}
	\right).
	\end{multline*}
	\normalsize
	The square root 
in the last line has to be correctly understood, as it will be explained 
in the remaining part of the proof.

	Since Cauchy transform $G$ of $\U = (\i+\X)^{-1/2}\,\Y\,(\i+\X)^{-1/2}$ satisfies  (see \eqref{for:G})
	$${G(z) = \frac{1}{z}\, A\left(\frac{1}{z}\right)},$$
	then 
	\bel{cauchy1}
		G(z) =   \frac{\bar\beta}{2 (-1+\bar\alpha\bar \beta) }
		\left(1 +\frac{\bar\alpha }{1+z}-\frac{ \bar\beta^{-1}+a
		-\bar\alpha}{z} +\frac{\bar\beta^{-1}}{z\,(1+z)}\sqrt{p_1(z)-p_2(z)}\right),
	\ee
	where 
	\begin{equation*}
		\begin{split}
		p_1(z) &=\left[-\bar\beta\,z^2 + z\, \bar\beta\,(\bar\beta^{-1} + a - 1-2\bar\alpha) + 1 + a\bar\beta -\bar\alpha\bar\beta\right]^2,\\
		p_2(z) &= 4\bar \beta\, (-1+\bar\alpha \bar\beta) \,z \,(1+z)\, (1-a+\bar\alpha+\alpha_1+z).
		\end{split}
	\end{equation*}
	We have to find an admissible set of parameters $\bar\alpha, \bar\beta, a, \alpha_1$, 
	such that $G$ is Cauchy transform of~the probabilistic measure associated with $\U$. 
	For that reason we analyse the polynomial $p(\cdot) := p_1(\cdot)-p_2(\cdot)$, which is under the square root in~formula \eqref{cauchy1}. 
	In the remainder of the proof, we will
	match formula \eqref{cauchy1} with the expression \eqref{cauchy} for the Cauchy transform of a free Kummer distribution,
	in particular identifying appropriate roots of the polynomials $p = p_1 - p_2$ as the boundary of
	the support of this distribution.
	It is known that $G$ is analytic on $\C^+$ and the image $G\left(\C^+\right)=\C^-$.
	We assumed that $\star$-distribution of $\U$ is supported in $(0, \infty)$. 
	These facts imply, that $p$ 
	\begin{itemize}
		\item[$(\star)$]  does not have roots in $\C^+$ and so it does not have complex roots at all;
		thus it has $4$ (possibly multiple) real roots;
		\item[$(\star\star)$] can not be negative for negative (real) arguments:  
		this follows from Stieltjes formula~\eqref{stielt}.
	\end{itemize}
	
	The roots of $p_2$ are: $z_1 = 0$, $z_2 = -1$ and $z_3 = a - \bar\alpha - \alpha_1 - 1$. 
	Since ${\fif(\V|\U) = \fif(\X+\Y - \U|\U)}$ and $\alpha_1 = \fif(\U)$ and~$a=\beta_0 =\fif(\Y)$,
	 then from \eqref{jed}: $$ a - \bar\alpha - \alpha_1 = -\fif(\X)<0.$$
	Therefore $z_3<-1$, which will be important later on in the proof. 
	Now we sketch the graph of $p_2$ -- see the first panel in Fig.~\eqref{w_p2}.
	\begin{center}\begin{figure}[t!]\centering{
		\includegraphics[scale=0.7]{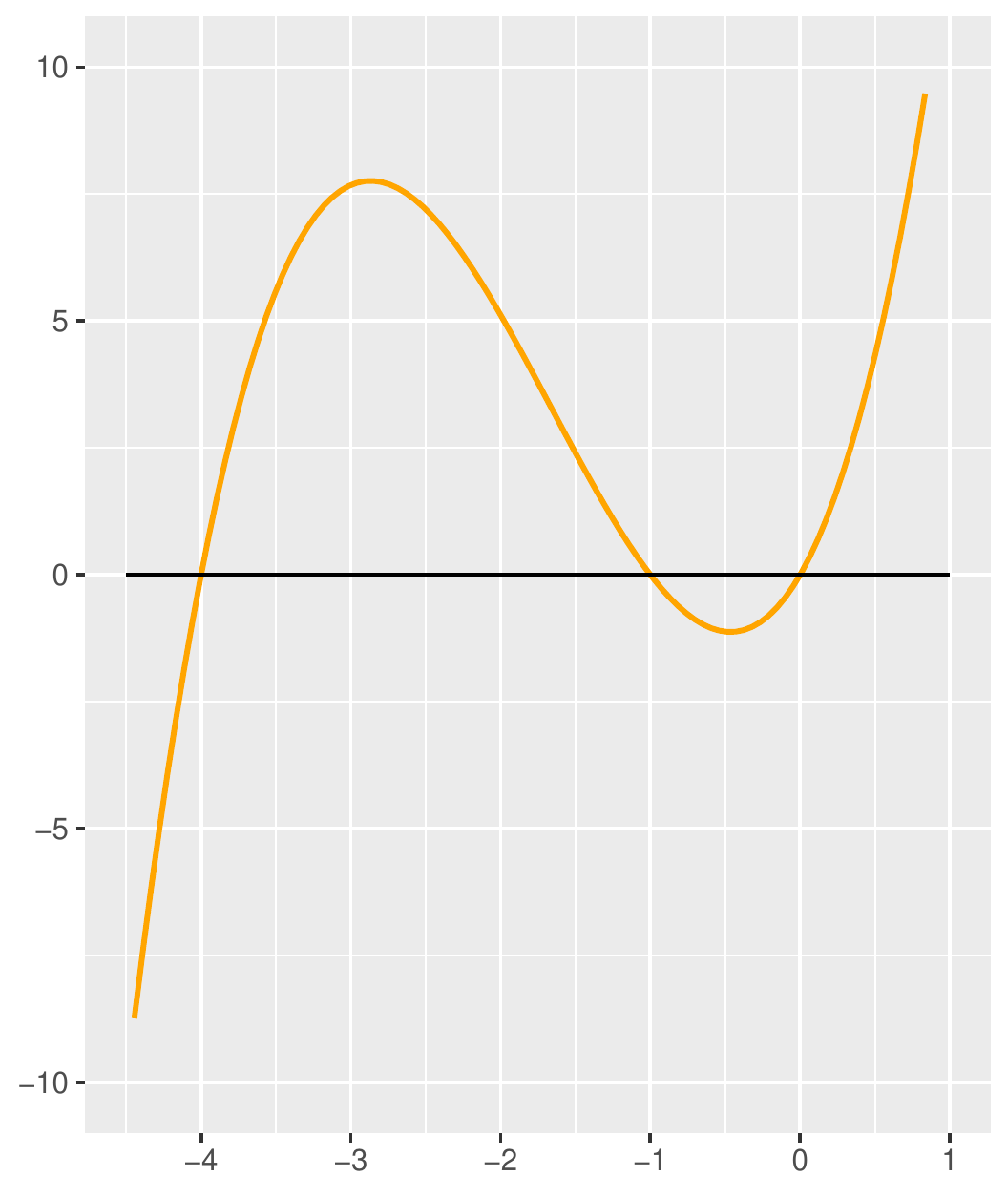}\hspace{14pt}
		\includegraphics[scale=0.7]{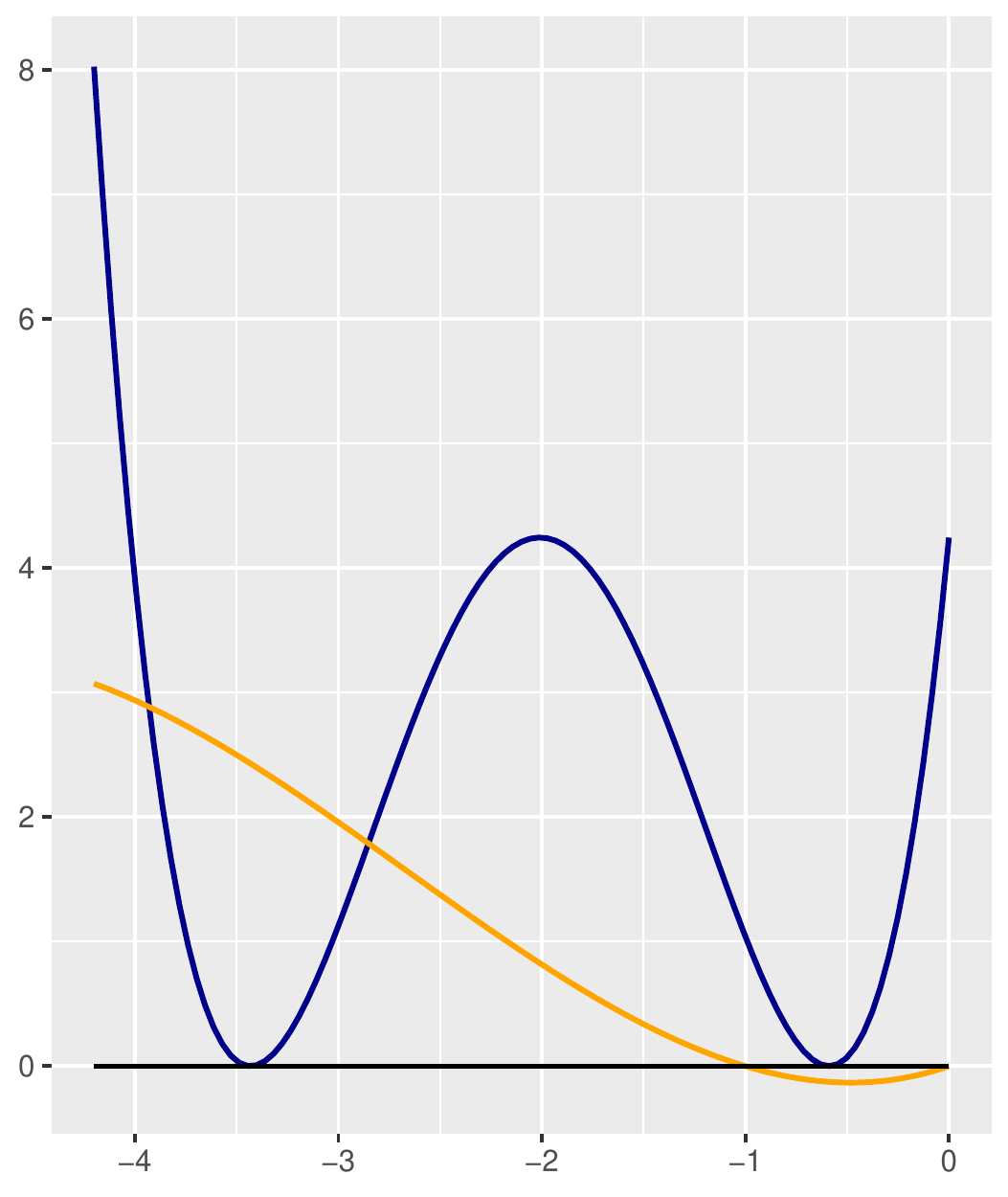}
		\caption{Left: $p_2$ (orange) for $\bar\alpha=4.5$, $\bar\beta=0.4$, $a=3$, $\alpha_1=1.5$. 
		Right: $p_1$ (blue)and $p_2$ (orange) for $\bar\alpha=1$, $\bar\beta=1.02$, $a=-2$, $\alpha_1=3$. \label{w_p2}}
	}
	\end{figure}\end{center}
	\vspace*{-0.9cm}
	We will show that $p_1$ has two double real roots now.
	Let
	$$f(z):= z^2 -z\, (d - 1-\bar\alpha) - d,$$
	where $d = \bar\beta^{-1} + a -\bar\alpha$.
	Then $p_1(z) = \bar\beta^2 f^2(z)$ and $f$ has the same roots as $p_1$.
	Note, that the discriminant of $f$ is positive. Indeed:
	\begin{equation*}
		\Delta_f  = (d - 1-\bar\alpha) ^2 + 4d
				   = (d+1-\bar\alpha)^2 + 4\bar\alpha>0,
	\end{equation*}
	where the last equality holds because $\bar\alpha>0$. 
	
	We denote roots of $p_1$: $\zeta_1<\zeta_2$. 
	If  $\zeta_1\cdot\zeta_2>0$, then Viete's formula implies
	$\bar\alpha - a -\bar\beta^{-1} >0$. 
	Then if $\zeta_1,\, \zeta_2>0$, we have $-(\bar\alpha - a  - \bar\beta^{-1}) - 1 -\bar\alpha=\zeta_1 + \zeta_2 >0$. 
	But this contradicts to the fact that $\bar\alpha>0$. 
	On the other hand, when $\zeta_1, \zeta_2<0$, then $p_1$ and $p_2$ do not intersect (which is contrary to $(\star)$) 
	or they intersect in the second quarter of the plane. 
	It contradicts  $(\star\star)$.
	This case is presented in the second panel of Fig.~\ref{w_p2}. 

	So we conclude that  ~$\zeta_1<0<\zeta_2$. 
	Or equivalently 
	\bel{pom2}
		\bar\beta^{-1} + a -\bar\alpha = -\zeta_1\cdot\zeta_2>0.
	\ee
	Now, we want to show that $p$ has two (distinct) positive roots and a negative double root.
	Suppose that $p$ has four different roots.
	Since $\zeta_1<0<\zeta_2$ and $p_2(x)>0$ for $x>0$,  
	then two roots of $p=p_1-p_2 $ are negative. 
	We denote them $y_1<y_2$. 
	So $p$ is negative in the interval $(y_1,\, y_2)$ (see Fig.~\ref{theugly}), which is contrary to $(\star\star)$. 
	This implies $p$ has double root $x_0$ 
	($p_1$ and $p_2$ are tangents at $x_0$). 
	Since $p_1$ and $p_2$ cannot be tangent outside interval $(z_3, -1)$, then $x_0<-1$. This is in the second chart in Fig.~\ref{theugly}.
	The other points of intersection of $p_1$ and $p_2$ are positive $x_1<x_2$. 	
	\begin{center}\begin{figure}[t!]\centering{
		\includegraphics[scale=0.7]{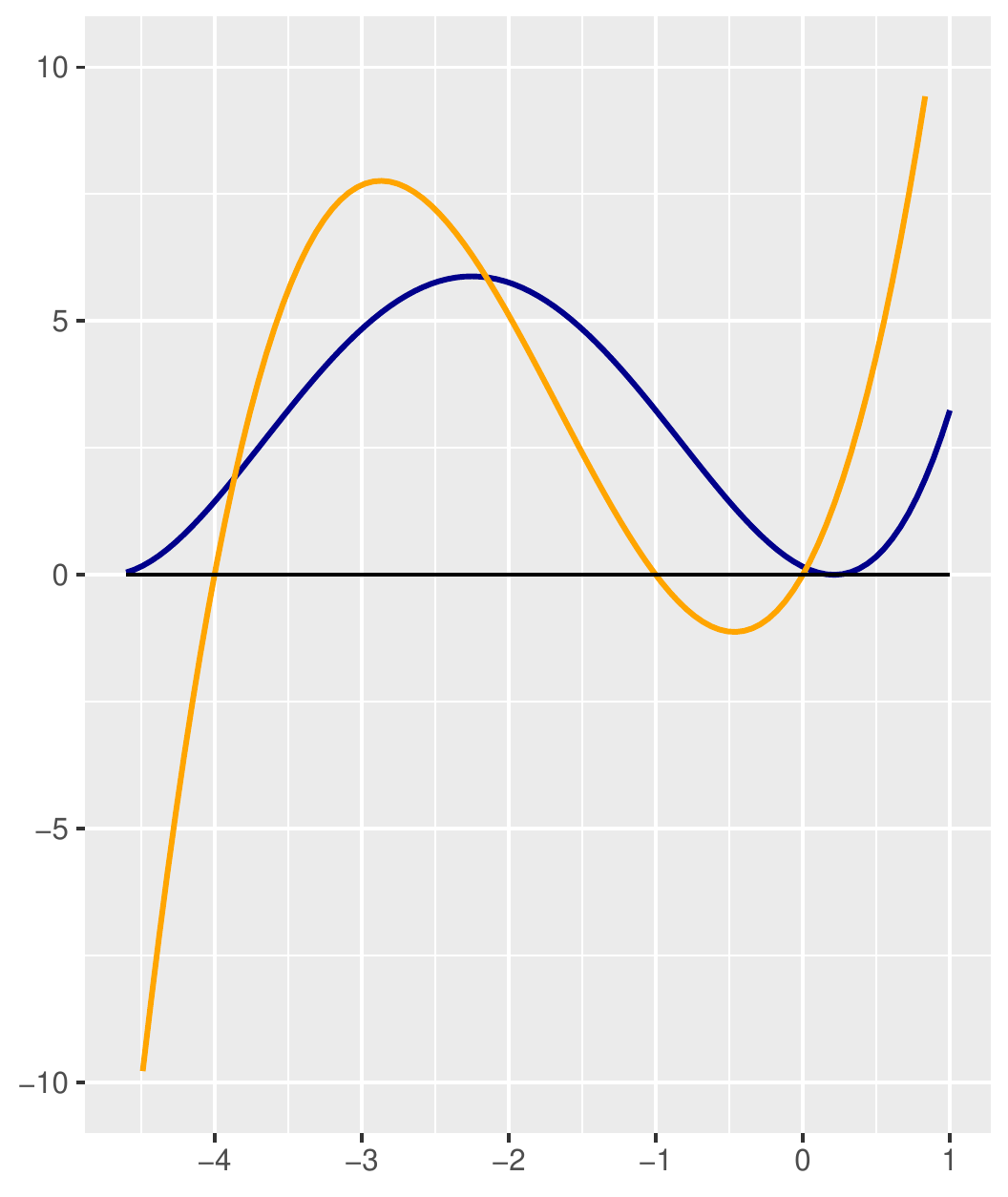}\hspace{8pt}
		\includegraphics[scale=0.7]{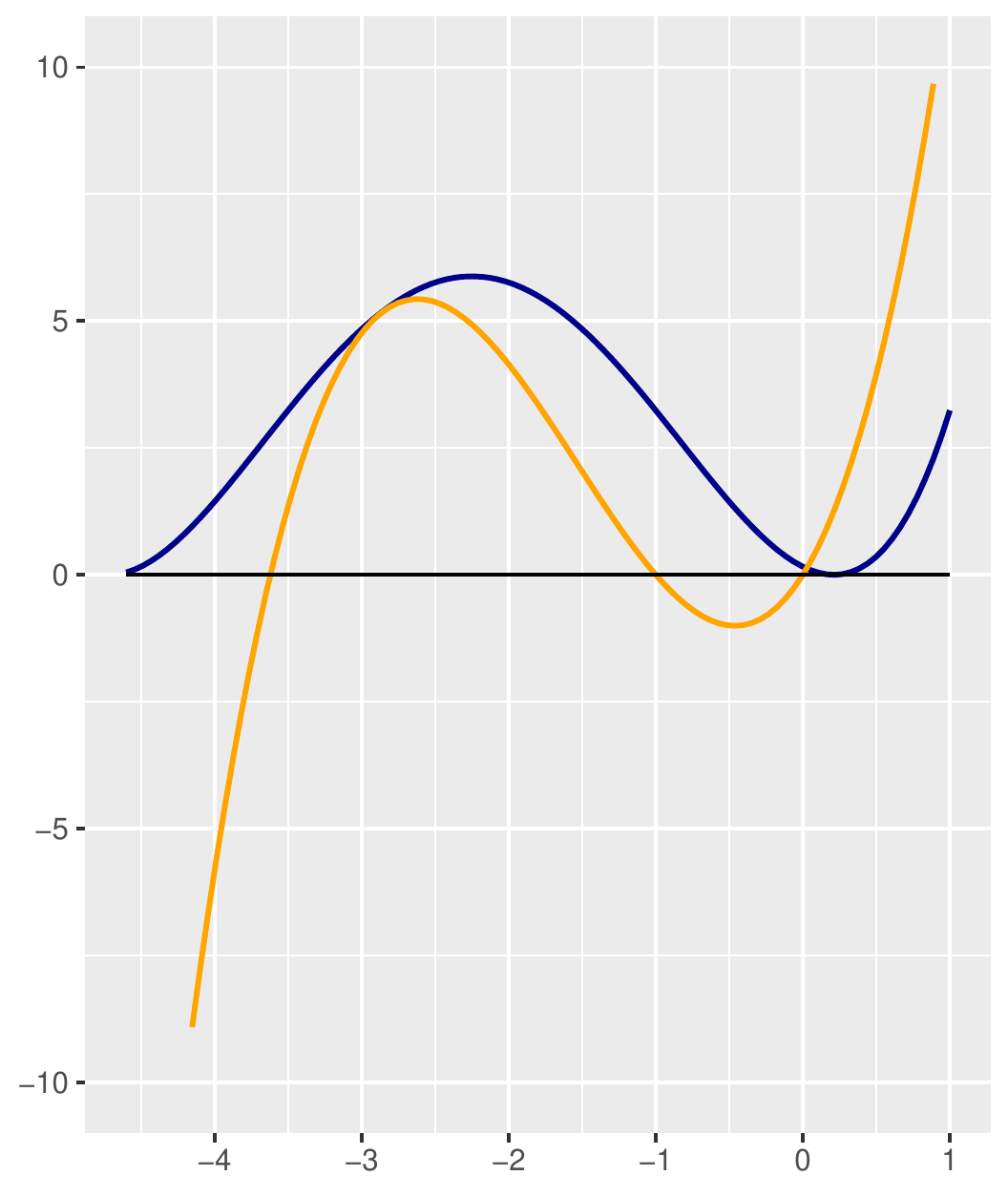}
		\caption{Polynomials $p_1$ (blue) and $p_2$ (orange) depending on parameters: 
		$\alpha_1 = 1.5$ on the left chart and $\alpha_1 = 1.12$ on the right chart 
		(other parameters are fixed: $\alpha = 4.5,\, \beta = 0.4,\, a=3$).}
		\label{theugly}}
	\end{figure}\end{center}
	\vspace*{-0.8cm}
	We have
	$$p(z)=(p_1-p_2)(z)= a_0+a_1\, z + a_2\, z^2 + a_3\, z^3 + a_4\, z^4,$$
	where
	\begin{equation*}
	\begin{split}
		a_0&=(1+a \bar\beta-\bar\alpha\bar\beta)^2,\\
		a_1&=2 \left[1+(1-\bar\alpha+2 \alpha_1)\, \bar\beta+\left(a^2-a \,(1+\bar\alpha)-\bar\alpha\, (1+2 \alpha_1)\right)\,\bar \beta^2\right],\\
		a_2& = 1+(4-2 a+4 \alpha_1)\, \bar\beta+\left(1-4 a+a^2-2\, \bar\alpha-4\, \bar\alpha\alpha_1\right) \,\bar\beta^2,\\
		a_3&= 2 \bar\beta\,(1+\bar\beta-a \bar\beta), \\
		a_4 &= \bar\beta^2.
	\end{split}
	\end{equation*}
	To find out how   $x_0,\, x_1$ and $x_2$ depend on $\bar\alpha$, $\bar\beta$ and $a$ 
	we use Viete's formula once again.
	We have
	\bel{pier1}
	x_1\, x_2\, x_0^2 = \frac{a_0}{a_4} = \left(\frac{1+a\bar\beta -\bar \alpha\bar\beta}{\bar\beta}\right)^2,
	\ee
	\bel{pier2}
	x_1+x_2 +2x_0 = -\frac{a_3}{a_4} = -2\,(1- a+\bar\beta^{-1}).
	\ee
	Since $x_0<-1<0$ and \eqref{pom2} holds, then from Eq. \eqref{pier1}
	\bel{gw_1}
	x_0 =\frac{- |\bar\beta^{-1} +a-\bar\alpha|}{\sqrt{x_1 x_2}} =-\frac{\bar\beta^{-1} +a-\bar\alpha}{\sqrt{x_1 x_2}}.
	\ee
	Let $q(z) := p(z-1)$. Then
	$$q(z) =  a_0' + a_1'\, z +a_2'\, z^2 + a_3'\, z^3 + a_4'\, z^4,$$
	where 
	\begin{equation*}
	\begin{split}
		a_0' &= \bar\alpha^2\bar\beta^2,\\
		a_1' &= -2\bar \beta \,\left[2 \alpha_1+a (-2+\bar\alpha\bar\beta)-\bar\alpha(-1+\bar\beta+2 \alpha_1\bar\beta)\right],\\
		a_2' &= 1-2\bar\beta\, (1+a-2 \alpha_1) +\left(1+2 a+a^2-2\bar \alpha-4 \bar\alpha\alpha_1\right) \bar\beta^2,\\
		a_3' &= -2\bar \beta\, (-1 + \bar\beta + a \bar\beta), \\
		a_4' &= \bar\beta^2.
	\end{split}
	\end{equation*}
	Notice that if $q(x')=0$, then $x = x'-1$ is a root of $p$. 
	So roots of $q$ are exactly $x_0+1$ (double root), $x_1+1$ and $x_2+1$. 
	Again from Viete's formula we have 
	\bel{pier3}
	(x_1+1)\,(x_2 +1)\, (x_0+1)^2 = \frac{a_0'}{a_4'} = \bar\alpha^2.
	\ee
	
	From Eq. \eqref{pier3} and again from $x_0<-1$, we have
	\bel{gw_2} x_0 +1=- \frac{\bar\alpha}{\sqrt{(x_1+1)(x_2+1)}}.\ee
	Let us recall that roots $x_1$ and $x_2$ are the boundary of the support of $\U=(\i+\X)^{-1/2}\,\Y\,(\i+\X)^{-1/2}$. 
	Combining \eqref{pier2} with~\eqref{gw_2}, we have 
	\bel{supp2}
		\frac{x_1+x_2}{2}  - a+\bar\beta^{-1} -  \frac{\bar\alpha}{\sqrt{(x_1+1)(x_2+1)}}  = 0.
	\ee
	From \eqref{gw_1} and~\eqref{gw_2} we have:
	\bel{supp3}\frac{\bar\alpha}{\sqrt{(x_1+1)(x_2+1)}} - \frac{\bar\beta +a - \bar\alpha}{\sqrt{x_1x_2}} + 1 = 0. \ee
	These are exactly conditions from \eqref{u2} for the boundary points of the support of free Kummer distribution with parameters 
	$\gamma = \bar\beta/(\bar\alpha\bar\beta - 1)$, $\alpha = (1/\bar\beta - \bar\alpha +a)\,\gamma+1 = a\gamma$ and $\beta = \bar\alpha\gamma$.
	So we have
	$$p(z) = (z-x_1)(z-x_2)(z-x_0)^2$$ 
	and from~\eqref{cauchy1}
	\begin{equation*}
	\begin{split}
	G(z) &= \frac{1}{2}\left[\gamma - \frac{\alpha-1}{z} +\frac{\beta}{1+z}+\sqrt{(z-x_1)(z-x_2)}\left(\frac{\gamma}{1+z} - \frac{\alpha-1}{\sqrt{x_1x_2} \,z\,(1+z)}\right)\right]
	\\
	&=  \frac{1}{2}\left[\gamma - \frac{\alpha-1}{z} +\frac{\beta}{1+z}+\sqrt{(z-x_1)(z-x_2)}\left(\frac{\beta}{(1+z)\sqrt{(x_1+1)(x_2+1)}} -\frac{\alpha-1}{z\,\sqrt{x_1x_2}}\right)\right],
	\end{split}
	\end{equation*}
	where $x_1$ and $x_2$ are such that \eqref{supp2} and \eqref{supp3} hold. 
	Since the support of $\U$ is bounded, we choose the main branch of square root
	(see the second part of the proof of Lem.~\ref{cauchykummer} for further reasoning).
	We finally obtained Cauchy transform of free Kummer with parameters $\alpha$, $\beta$, $\gamma$.

	\textbf{Step 3.}
	Having the distributions of $\Y$ and $\U$ already identified, we can calculate $\X$ by its $S$-transform. 
	Let $S_\mathbb{Z}$ denote the $S$-transform of random variable $\mathbb{Z}$. 	
	Given that $\X$ and $\Y$ are free, 
	Thm.~\ref{T:S} implies
	$$S_\Y(z)\, S_{(\i+\X)^{-1}}(z) = S_\U(z)$$
	and this equality uniquely determines distribution of $\X$. 
	Theorem \ref{wprost} implies that
	${\X\sim f\K (\beta, \alpha, \gamma)}$.
\end{proof}

\section*{Acknowledgments}
This research was supported by the grant \texttt{2016/21/B/ST1/00005} of National Science Center, Poland.

\end{document}